# KARDAR-PARISI-ZHANG UNIVERSALITY

IVAN CORWIN

This is the article by the same name which was published in the March 2016 issue of the AMS Notices. Due to space constraints, some references contained here were not provided in the published version were not provided. The figures in this article were made collaboratively with William Casselman.

## 1. Universality in random systems

Universality in complex random systems is a striking concept which has played a central role in the direction of research within probability, mathematical physics and statistical mechanics. In this article we will describe how a variety of physical systems and mathematical models, including randomly growing interfaces, certain stochastic PDEs, traffic models, paths in random environments, and random matrices all demonstrate the same universal statistical behaviors in their long-time / large-scale limit. These systems are said to lie in the *Kardar-Parisi-Zhang (KPZ) universality class*. Proof of universality within these classes of systems (except for random matrices) has remained mostly elusive. Extensive computer simulations, non-rigorous physical arguments / heuristics, some laboratory experiments and limited mathematically rigorous results provide important evidence for this belief.

The last fifteen years have seen a number of breakthroughs in the discovering and analysis of a handful of special *integrable probability systems* which, due to enhanced algebraic structure, admit many exact computations and ultimately asymptotic analysis revealing the purportedly universal properties of the KPZ class. The structures present in these systems generally originate in representation theory (e.g. symmetric functions), quantum integrable systems (e.g. Bethe ansatz), algebraic combinatorics (e.g. RSK correspondence) and the techniques in their asymptotic analysis generally involves Laplaces method, Fredholm determinants, or Riemann-Hilbert problem asymptotics.

This article will focus on the phenomena associated with the KPZ universality class [14] and highlight how certain integrable examples expand the scope of and refine the notion of universality. We start by providing a brief introduction to the Gaussian universality class and the integrable probabilistic example of random coin flipping and the random deposition model. A small perturbation to the random deposition model leads us to the ballistic deposition model and the KPZ universality class. The ballistic deposition model fails to be integrable, thus to gain an understanding of its long-time behavior and that of the entire KPZ class, we turn to the corner growth model. The rest of the article focuses on various sides of this rich model – its role as a random growth process, its relation to the KPZ stochastic PDE, its interpretation in terms of interacting particle systems, and its relation to optimization problems involving paths in random environments. Along the way, we include some other generalizations of this process whose integrability springs from the same sources. We close the article by reflecting upon some open problems.

A survey of the KPZ universality class and all of the associated phenomena and methods developed or utilized in its study is far too vast to be provided here. This article presents only one of many stories and perspectives regarding this rich area of study. To even provide a representative cross-section of references is beyond this scope. Additionally, though we will discuss integrable examples, we will not described the algebraic structures and methods of asymptotic analysis behind them (despite their obvious importance and interest). Some recent references which review some these structures include [8, 11, 12, 15, 7, 26] and references therein. On the more physics oriented side, the collection of reviews and books [3, 14, 20, 23, 24, 26, 31, 32] provides some idea of the scope of the study of the KPZ universality class and the diverse areas upon which it touches.

We start now by providing an overview of the general notion of universality in the context of the simplest and historically first example – fair coin flipping and the Gaussian universality class.

## 2. Gaussian universality class

Flip a fair coin $N$ times. Each string of outcomes (e.g. head, tail, tail, tail, head) has equal probability $2^{-N}$. Call $H$ the (random) number of heads and let $\mathbb{P}$ denote the probability distribution for this sequence of coin flips. Counting shows that

$$\mathbb{P}(H = n) = 2^{-N}\binom{N}{n}.$$

Since each flip is independent, the expected number of heads is $N/2$. Bernoulli (1713) proved that $H/N$ converges to $1/2$ as $N$ goes to infinity. This was the first example of a *law of large numbers*. Of course, this does not mean that if you flip the coin 1000 times, you will see exactly 500 heads. Indeed, in $N$ coin flips one expects the number of heads to vary randomly around the value $N/2$ in the scale $\sqrt{N}$. Moreover, for all $x \in \mathbb{R}$,

$$\lim_{N\to\infty} \mathbb{P}\Big(H < \tfrac{1}{2}N + \tfrac{1}{2}\sqrt{N}\,x\Big) = \int_{-\infty}^{x} \frac{e^{-y^2/2}}{\sqrt{2\pi}} dy.$$

De Moivre (1738), Gauss (1809), Adrain (1809), and Laplace (1812) all participated in the proof of this result. The limiting distribution is known as the Gaussian (or sometimes normal or bell curve) distribution.

A proof of this follows from asymptotics of $n!$, as derived by de Moivre (1721) and named after Stirling (1729). Write

$$n! = \Gamma(n+1) = \int_0^\infty e^{-t} t^n dt = n^{n+1} \int_0^\infty e^{nf(z)} dz$$

where $f(z) = \log z - z$ and the last equality is from the change of variables $t = nz$. The integral is dominated, as $n$ grows, by the maximal value of $f(z)$ on the interval $[0, \infty)$. This occurs at $z = 1$, thus expanding $f(z) \approx -1 - \frac{(z-1)^2}{2}$ and plugging this into the integral yields the final expansion

$$n! \approx n^{n+1} e^{-n} \sqrt{2\pi/n}.$$

This general route of writing exact formulas for probabilities in terms of integrals and then performing asymptotics is quite common to the analysis of integrable models in the KPZ universality class – though those formulas and analyses are considerably more involved.

The universality of the Gaussian distribution was not broadly demonstrated until work of Chebyshev, Markov and Lyapunov around 1900. The *central limit theorem* showed that the exact nature of coin flipping is immaterial – any sum of independent identically distributed (iid) random variables with finite mean and variance will demonstrate the same limiting behavior.

**Theorem 2.1.** *Let $X_1, X_2, \ldots$ be iid random variables of finite mean $m$ and variance $v$. Then for all $x \in \mathbb{R}$,*

$$\lim_{N\to\infty} \mathbb{P}\bigg(\sum_{i=1}^N X_i < mN + v\sqrt{N}x\bigg) = \int_{-\infty}^{x} \frac{e^{-y^2/2}}{\sqrt{2\pi}} dy.$$

Proofs of this result use different tools than the exact analysis of coin flipping and much of probability theory deals with the study of Gaussian processes which arise through various generalizations of the CLT. The Gaussian distribution is ubiquitous and, as it is the basis for much of classical statistics and thermodynamics, it has had immense societal impact.

## 3. Random versus ballistic deposition

The *random deposition model* is one of the simplest (and least realistic) models for a randomly growing one-dimensional interface. Unit blocks fall independently and in parallel from the sky above each site of $\mathbb{Z}$ according to exponentially distributed waiting times (see Figure 1). Recall that a random variable $X$ has exponential distribution of rate $\lambda > 0$ (or mean $1/\lambda$) if $\mathbb{P}(X > x) = e^{-\lambda x}$. Such random variables are characterized by the memoryless property – conditioned on the event that $X > x$, $X - x$ still has the exponential distribution of the same rate. Consequently, the random deposition model is *Markov* – its future evolution only depends on the present state (and not on its history).

The random deposition model is quite simple to analyze since each column grows independently. Let $h(t, x)$ record the height above site $x$ at time $t$ and assume $h(0, x) \equiv 0$. Define random waiting times $w_{x,i}$ to be the time for the $i$-th block in column $x$ to fall. For any $n$, the event $h(t, x) < n$ is equivalent to $\sum_{i=1}^n w_{x,i} > t$. Since the $w_{x,i}$ are iid,

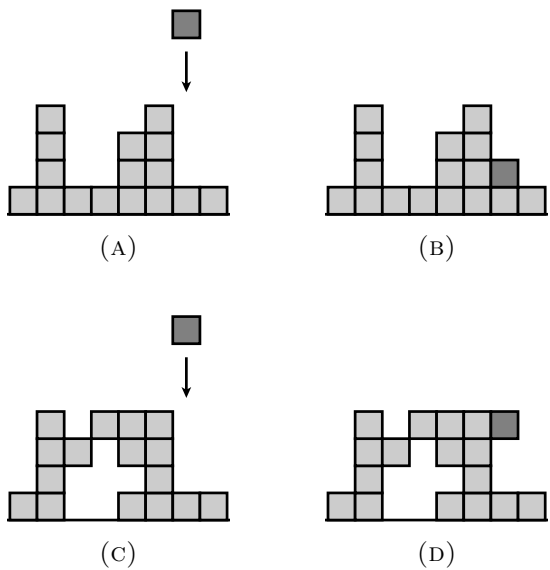

Figure 1. (A) and (B) illustrate the *random deposition model* and (C) and (D) illustrate the *ballistic deposition model*. In both cases, blocks fall from above each site with independent exponentially distributed waiting times. In the first model, they land at the top of each column whereas in the second model they stick to the first edge to which they become incident.

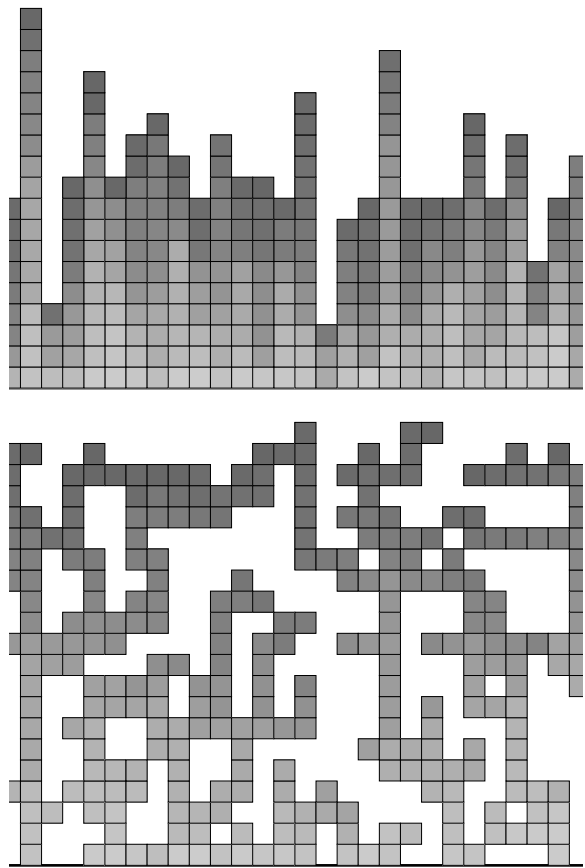

Figure 2. Simulation of random (top) versus ballistic (bottom) deposition models driven by the same process of falling blocks. The ballistic model grows much faster, and has a smoother more spatially correlated top interface.

the law of large numbers and central limit theory apply here. Assuming $\lambda = 1$,

$$\lim_{t\to\infty} \frac{h(t,x)}{t} = 1, \quad \text{and} \quad \lim_{t\to\infty} \frac{h(t,x) - t}{t^{1/2}} \Rightarrow N(x)$$

jointly over $x \in \mathbb{Z}$, where $\{N(x)\}_{x\in\mathbb{Z}}$ is a collection of iid standard Gaussian random variables. The top of Figure 2 shows a simulation of the random deposition model. The linear growth speed and lack of spatial correlation are quite evident. The fluctuation of this model are said to be in the Gaussian universality class since they grow like $t^{1/2}$, with Gaussian limit law and trivial transversal correlation length scale $t^0$. In general, fluctuation and transversal correlation exponents, as well as limiting distributions constitutes the description of a universality class and all models which match these limiting behaviors are said to lie in the same universality class.

While the Gaussian behavior of this model is resilient against changes in the distribution of the $w_{x,i}$ (owing to the CLT), generic changes in the nature of the growth rules shatter the Gaussian behavior. The *ballistic deposition (or sticky block) model* was introduced by Vold [38] in 1959 and, as one expects in real growing interfaces, displays spatial correlation. As before, blocks fall according to iid exponential waiting times, however, now a block will stick to the first edge against which it becomes incident. This mechanism is illustrated in Figure 1. This creates overhangs and we define the height function $h(t,x)$ as the maximal height above $x$ which is occupied by a box. How does this microscopic change manifest itself over time?

It turns out that sticky blocks radically changes the limiting behavior of this growth process. The bottom of Figure 2 records one simulation of the process. Seppäläinen [29] gave a proof that there is still an overall linearly growth rate. Moreover by considering a lower bound by a width two system, one can see that this velocity exceeds that of the random deposition model. The exact value of this rate, however, remains unknown.

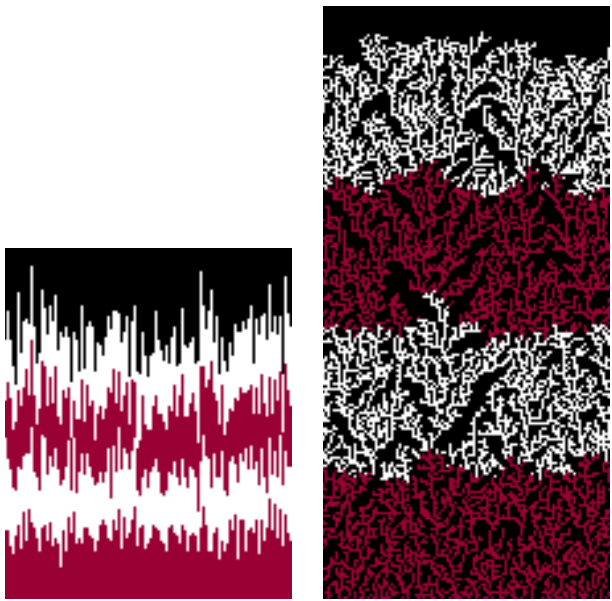

FIGURE 3. Simulation of random (left) versus ballistic (right) deposition models driven by the same process of falling blocks and run for a long time. The red and white colors represent different epochs of time in the simulation. The size of boxes in both figures are the same.

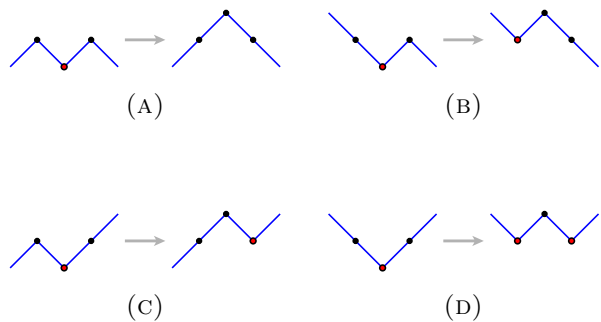

FIGURE 4. Various possible ways that a local minimum can grow into a local maximum. The red dot represent the local minimum at which growth may occur.

The simulation in Figure 2 (as well as the longer time results displayed in Figure 3) also shows that the scale of fluctuations of $h(t, x)$ is smaller than in random deposition, and that the height function remains correlated transversally over a long distance. There are exact conjectures for these fluctuations. They are supposed to grow like $t^{1/3}$ and demonstrate a non-trivial correlation structure in a transversal scale of $t^{2/3}$. Additionally, precise predictions exist for the limiting distributions. Up to certain (presently undetermined) constants $c_1, c_2$, the sequence of scaled heights $c_2 t^{-1/3}(h(t, 0) - c_1 t)$ should converge to the so-called *Gaussian orthogonal ensemble (GOE) Tracy-Widom* distributed random variable. The Tracy-Widom distributions can be thought of as modern-day bell curves, and their names GOE or GUE (for *Gaussian unitary ensemble*) come from the random matrix ensembles in which these distributions were first observed by Tracy-Widom [33, 34].

Ballistic deposition does not seem to be an integrable probabilistic system, so where do these precise conjectures come from? The exact predictions come from the analysis of a few similar growth processes which just happen to be integrable! Ballistic deposition shares certain features with these models which are believed to be key for membership in the KPZ class:

- Locality: Height function change depends only on neighboring heights.
- Smoothing: Large valleys are quickly filled.
- Non-linear slope dependence: Vertical effective growth rate depends non-linearly on local slope.
- Space-time independent noise: Growth is drive by noise which quickly decorrelates in space / time and is not heavy tailed.

It should be made clear that a proof of the KPZ class behavior for the ballistic deposition model is far beyond what can be done mathematically (though simulations strongly suggest that the above conjecture is true).

## 4. Corner growth model

We come to the first example of an integrable probabilistic system in the KPZ universality class – the *corner growth model*. The randomly growing interface is modeled by a height function $h(t, x)$ which is continuous, piece-wise linear, and composed of $\sqrt{2}$-length line increments of slope $+1$ or $-1$, changing value at integer $x$. The height function evolves according to the Markovian dynamics that each local minimum of $h$ (looking like ∨) turns into a local maximum (looking like ∧) according to an exponentially distributed waiting time. This happens independently for each minimum. This change in height function can also be thought of as adding boxes (rotated by 45°). See Figures 4 and 5 for further illustration of this model.

*Wedge* initial data means that $h(0,x) = |x|$ while *flat* initial data (as considered for ballistic deposition) means that $h(0,x)$ is given by a periodic sawtooth function which goes between height 0 and 1. We will focus on wedge initial data. In 1980, Rost [27] proved a law of large numbers for the growing interface when time, space and the height function are scaled by the same large parameter $L$.

**Theorem 4.1.** *For wedge initial data,*

$$\lim_{L\to\infty} \frac{h(Lt, Lx)}{L} = \mathfrak{h}(t,x) := \begin{cases} t\frac{1-(x/t)^2}{2} & |x| < t, \\ |x| & |x| \geq t. \end{cases}$$

Figure 6 displays the result of a computer simulation wherein the limiting parabolic shape is evident. The function $\mathfrak{h}$ is the unique viscosity solution to the Hamilton-Jacobi equation

$$\frac{\partial}{\partial t}\mathfrak{h}(t,x) = \frac{1}{2}\Big(1 - \big(\frac{\partial}{\partial x}\mathfrak{h}(t,x)\big)^2\Big).$$

This equation actually governs the evolution of the law of large numbers from arbitrary initial data.

The fluctuations of this model around the law of large numbers are what is believed to be universal. Figure 6 shows that the interface (blue) fluctuates around its limiting shape (red) on a fairly small scale, with transversal correlation on a larger scale. For $\epsilon > 0$, define the scaled and centered height function

$$h_\epsilon(t,x) := \epsilon^b h(\epsilon^{-z} t, \epsilon^{-1} x) - \frac{\epsilon^{-1} t}{2}$$

where the *dynamic scaling exponent* $z = 3/2$ and the *fluctuation exponent* $b = 1/2$. These exponents are easily remembered since they correspond with scaling time : space : fluctuations like 3 : 2 : 1. These are the characteristic exponents for the KPZ universality class. In 1999, Johansson [21] proved that for fixed $t$, as $\epsilon \to 0$, the random variable $h_\epsilon(t,0)$ converges to a GUE Tracy-Widom distributed random variable (see Figure 7). The work on the corner growth model was preceded by the work of Baik-Deift-Johansson [2] on the related model of the longest increasing subsequence in a random permutation. For that related model, two years later Prähofer-Spohn [25] computed the analog to the joint distribution of $h_\epsilon(t, x)$ for fixed $t$ and varying $x$. It was in [2] that the GUE Tracy-Widom distribution first appeared in combinatorial problems.

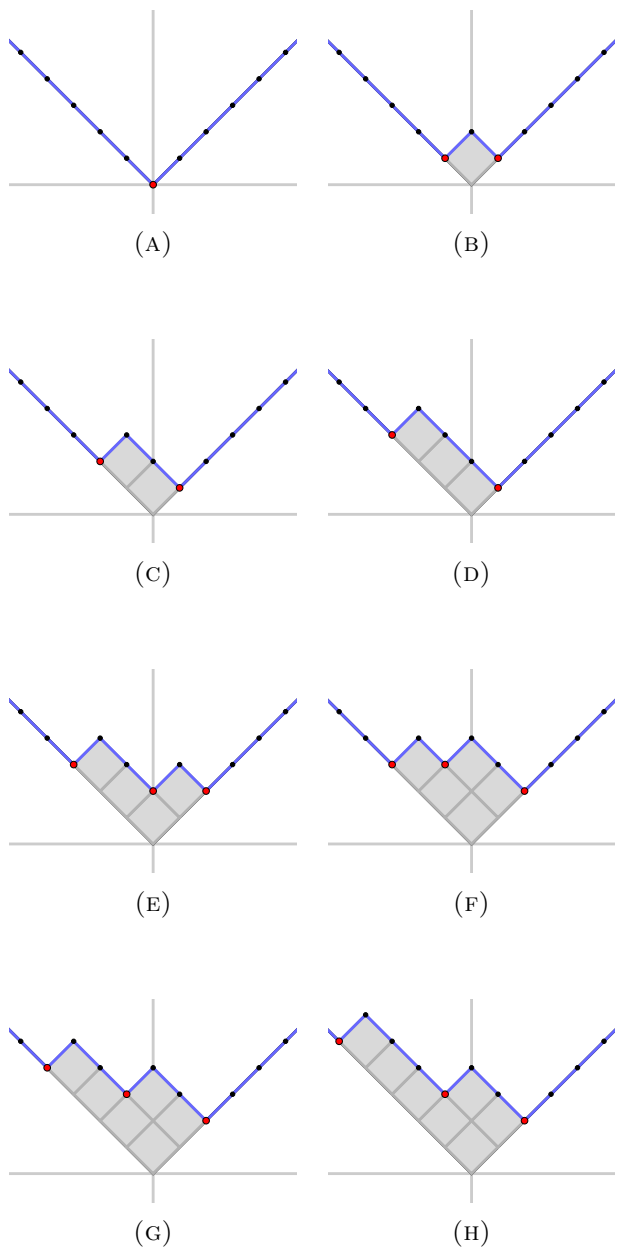

FIGURE 5. The corner growth model starts with an empty corner, as in (A). There is only one local minimum (the red dot) and after an exponentially distributed waiting time, this turns into a local maximum by filling in the site above it with a block, as in (B). In (B) there are now two possible locations for growth (the two red dots). Each one has an exponentially distributed waiting time. (C) corresponds to the case when the left local minimum grows before the right one. By the memoryless property of exponential random variables, once in state (C), we can think of choosing new exponentially distributed waiting times for the possible growth destinations. Continuing in a similar manner, we arrive at the evolution in (D) through (H).

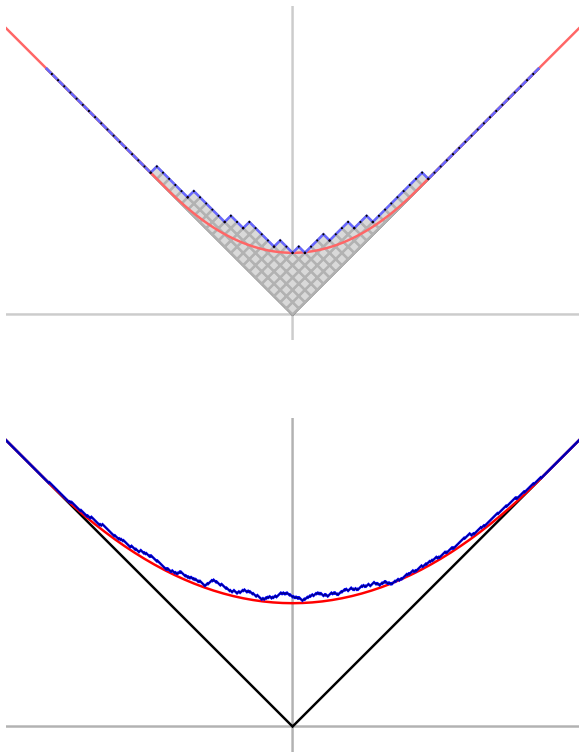

Figure 6. Simulation of the corner growth model. The top shows the model after a medium amount of time and the bottom shows it after a longer amount of time. The blue interface is the simulation while the red curve is the limiting parabolic shape. The blue curve has vertical fluctuations of order $t^{1/3}$ and decorrelates spatially on distances of order $t^{2/3}$.

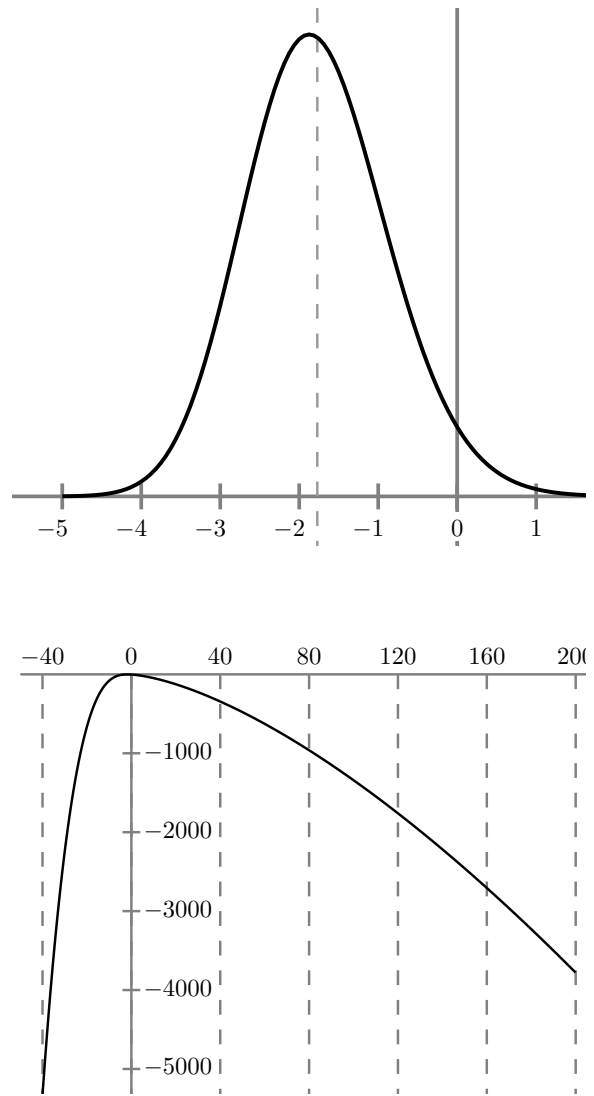

Figure 7. The density (top) and log of the density (bottom) of the GUE Tracy-Widom distribution. Though the density appears to look like a bell-curve (or Gaussian), this comparison is misleading. The mean and variance of the distribution are approximately $-1.77$ and $0.81$. The tails of the density (as shown in terms of the log of the density in the bottom plot) decay like $e^{-c_-|x|^3}$ for $x \ll 0$ and like $e^{-c_+ x^{3/2}}$ for $x \gg 0$, for certain positive constants $c_-$ and $c_+$. The Gaussian density decays like $e^{-cx^2}$ in both tails, with the constant $c$ related to the variance.

The entire scaled growth process $h_\epsilon(\cdot, \cdot)$ should have a limit as $\epsilon \to 0$ which would necessarily be a fixed point under the 3:2:1 scaling. The existence of this limit (often called the *KPZ fixed point*) remains conjectural. Still, much is known about the properties this limit should enjoy. It should be a stochastic process whose evolution depends on the limit of the initial data under the same scaling. The one-point distribution for general initial data, the multi-point and multi-time distribution for wedge initial data, and various aspects of its continuity are all understood. Besides the existence of this limit, what is missing is a useful characterization of the KPZ fixed point. Since the KPZ fixed point is believed to be the universal scaling limit of all models in the KPZ universality class and since corner growth enjoys the same key properties as ballistic deposition, one also is led to conjecture that ballistic deposition scales to the same fixed point and hence enjoys the same scalings and limiting distributions. The reason why the GOE Tracy-Widom distribution came up in our earlier discussion is that we were dealing with flat rather than wedge initial data.

One test of the universality belief is to introduce partial asymmetry into the corner growth model. Now we change local minimum into local maximum at rate $p$, and turn local maximum into local minimum at rate $q$ (all waiting times are independent and exponentially distributed, and $p+q=1$). See 8 for an illustration of this partially asymmetric corner growth model. Tracy-Widom [35, 36, 37] showed that so long as $p > q$, the same law of large numbers and fluctuation limit theorem holds for the partially asymmetric model, provided that $t$ is replaced by $t/(p-q)$. Since $p-q$ represents the growth drift, one simply has to speed up to compensate for this drift being smaller.

Clearly for $p \leq q$ something different must occur than for $p > q$. For $p = q$ the law of large numbers and fluctuations change nature. The scaling of time : space : fluctuations becomes 4 : 2 : 1 and the limiting process under these scalings becomes the stochastic heat equation with additive white-noise. This is the *Edwards-Wilkinson* (EW) universality class which is described by the stochastic heat equation with additive noise. For $p < q$ the process approaches a stationary distribution where the probability of having $k$ boxes added to the empty wedge is proportional to $(p/q)^k$.

So, we have observed that for any positive asymmetry the growth model lies in the KPZ universality class while for zero asymmetry it lies in the EW universality class. It is natural to wonder whether critically scaling parameters (i.e. $p - q \to 0$) one might encounter a crossover regime between these two universality classes. Indeed, this is the case, and the crossover is achieved by the KPZ equation which we now discuss.

## 5. KPZ EQUATION

The *KPZ equation* is written as

$$\frac{\partial h}{\partial t}(t,x) = \nu \frac{\partial^2 h}{\partial x^2}(t,x) + \tfrac{1}{2}\lambda\Big(\frac{\partial h}{\partial x}(t,x)\Big)^2 + \sqrt{D}\xi(t,x)$$

where $\xi(t,x)$ is Gaussian space-time white noise, $\lambda, \nu \in \mathbb{R}$, $D > 0$ and $h(t,x)$ is a continuous function of time $t \in \mathbb{R}_+$ and space $x \in \mathbb{R}$, taking values in $\mathbb{R}$. Due to the white-noise, one expects $x \mapsto h(t,x)$ to be only as regular in as Brownian motion. Hence, the non-linearity does not a priori make any sense (the derivative of Brownian motion has negative Hölder regularity). Bertini-Cancrini [5] provided the physically relevant notion of solution (called the *Hopf-Cole solution*) and showed how it arises from regularizing the noise, solving the (now well-posed) equation and then removing the noise and subtracting a divergence.

The equation contains the four key features mentioned earlier – the growth is local, depending on the Laplacian (smoothing), the square of the gradient (non-linear slope dependent growth), and white-noise (space-time uncorrelated noise). Kardar, Parisi, and Zhang [22] introduced their eponymous equation and 3 : 2 : 1 scaling prediction in 1986 in an attempt to understand the scaling behaviors of random interface growth.

How might one see the 3 : 2 : 1 scaling from the KPZ equation? Define $h_\epsilon(t,x) = \epsilon^b h(\epsilon^{-z}t, \epsilon^{-1}x)$, then $h_\epsilon$ satisfies the KPZ equation with scaled coefficients $\epsilon^{2-z}\nu$, $\epsilon^{2-z-b}\tfrac{1}{2}\lambda$ and $\epsilon^{b-\frac{z}{2}+\frac{1}{2}}\sqrt{D}$. It turns out that two-sided Brownian motion is stationary for the KPZ equation, hence any non-trivial scaling must respect the Brownian scaling of the initial data and thus have $b = 1/2$. Plugging this in, the only way to have no coefficient blow up to infinity, and not every term shrink to zero (as $\epsilon \to 0$) is to choose $z = 3/2$. This suggests that the plausibility of the 3 : 2 : 1 scaling. While this heuristic gives the right scaling, it does not provide for the scaling limit. The limit as $\epsilon \to 0$ of the equation (the inviscid Burgers equation where only the non-linearity survives) certainly does not govern the limit of the solutions. It remains something of a mystery as to exactly how to describe this limiting KPZ fixed point. The above heuristic says nothing of the limiting distribution of the solution to the KPZ equation, and there does not presently exist a simple way to see what this should be.

It took just under 25 years until Amir-Corwin-Quastel [1] rigorously proved that the KPZ equation is in the KPZ universality class. That work also computed an exact formula for the probability distribution of the solution to the KPZ equation – marking the first instance of a non-linear stochastic PDE for which this was accomplished. Tracy-Widom's work on the partially asymmetric corner growth model, and work of Bertini-Giacomin [6] which relates that model to the KPZ equation were the two main inputs in this development. See [14] for further details regarding this as well as the simultaneous exact but non-rigorous steepest descent work of Sasamoto-Spohn [28], and non-rigorous replica approach work of Calabrese-Le Doussal-Rosso [13], and Dotsenko [17].

The proof that the KPZ equation is in the KPZ universality class was part of an ongoing flurry of activity surrounding the KPZ universality class from a number of directions such as integrable probability [7, 15], experimental physics [32] and stochastic PDEs. For instance, Bertini-Cancrini's Hopf-Cole solution relies upon a trick (the Hopf-Cole transform) which linearizes the KPZ equation. Hairer [18], who had been developing methods to make sense of classically ill-posed stochastic PDEs, focused on the KPZ equation and developed a direct notion of solution which agreed with the Hopf-Cole one but did not require use of the Hopf-Cole transform trick. Still, this does not say anything about the distribution of solutions or their long-time scaling behaviors. Hairer's KPZ work set the stage for his development of regularity structures [19] – an approach to construction solutions of certain types of ill-posed stochastic PDEs – work for which he was awarded a Fields medal.

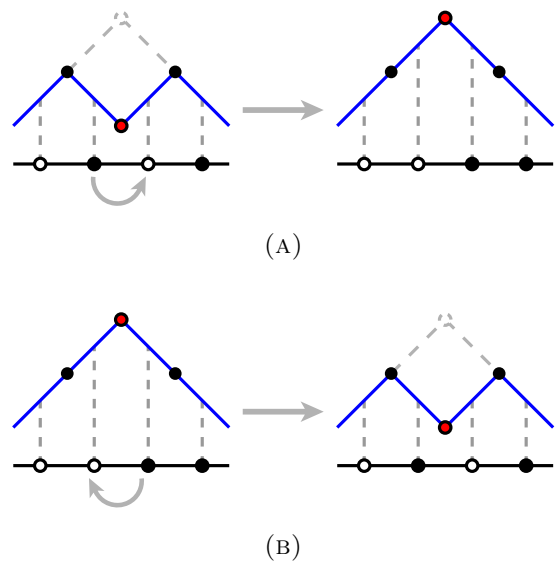

FIGURE 8. Mapping the partially asymmetric corner growth model to the partially asymmetric simple exclusion process. In (A), the red dot is a local minimum and it grows into a maximum. In terms of the particle process beneath it, the minimum corresponds to a particle followed by a hole, and the growth corresponds to said particle jumping into the hole to its right. In (B), the opposite is shown. The red dot is a local maximum and shrinks into a minimum. Correspondingly, there is a hole follows by a particle, and the shrinking results in the particle moving into the hole to its left.

## 6. Interacting particle systems

There is a direct mapping (see Figure 8) between the partially asymmetric corner growth model and the *partially asymmetric simple exclusion process* (generally abbreviated ASEP). Associate to every $-1$ slope line increment a particle on the site of $\mathbb{Z}$ above which the increment sits, and to every $+1$ slope line increment associate an empty site. The height function then maps onto a configuration of particles and holes on $\mathbb{Z}$, with at most one particle per site. When a minimum of the height function becomes a maxima, it corresponds to a particle jumping right by one into an empty site, and likewise when a maximum becomes a minimum, a particle jumps left by one into an empty site. Wedge initial data for corner growth corresponds with having all sites to the left of the origin initially occupied and all to the right empty – this is often called *step initial data* due to the step function in terms of particle density. ASEP was introduced in biology literature in 1968 by MacDonald-Gibbs-Pipkin as a model for RNA's movement during transcription. Soon after it was independently introduced within the probability literature in 1970 by Spitzer.

The earlier quoted results regarding corner growth immediately imply that the number of particles to cross the origin after a long time $t$ demonstrates KPZ class fluctuation behavior. KPZ universality would have that generic changes to this model should not change the KPZ class fluctuations. Unfortunately, such generic changes destroy the model's integrable structure. There are a few integrable generalizations discovered in the past five years which demonstrate some of the resilience of the KPZ universality class against perturbations.

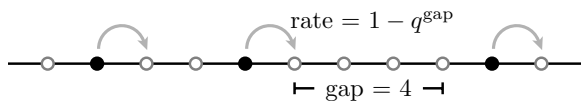

FIGURE 9. The $q$-TASEP, whereby each particle jumps one to the right after an exponentially distributed waiting time with rate given by $1-q^{\mathrm{gap}}$.

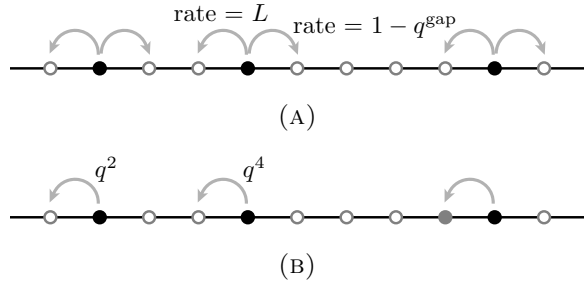

FIGURE 10. The $q$-pushASEP. As shown in (A), particles jump right according to the $q$-TASEP rates and left according to independent exponentially distributed waiting times of rate $L$. When a left jump occurs, it may trigger a cascade of left jumps. As shown in (B), the right-most particle has just jumped left by one. The next particle (to its left) instantaneously also jumps left by one with probability given by $q^{\mathrm{gap}}$ where gap is the number of empty sites between the two particles before the left jumps occurred (in this case gap = 4). If that next left jump is realized, the cascade continues to the next-left particle according to the same rule, otherwise it stops and no other particles jump left in that instant of time.

TASEP (the totally asymmetric version of ASEP) is a very basic model for traffic on a one-lane road in which cars (particles) move forward after exponential rate one waiting times, provided the site is unoccupied. A more realistic model would account for the fact that cars slow down as they approach the one in front. The model of $q$-TASEP does just that (Figure 10). Particles jump right according to independent exponential waiting times of rate $1-q^{\mathrm{gap}}$ where gap is the number of empty spaces to the next particle to the right. Here $q \in [0,1)$ is a different parameter than in the ASEP, though when $q$ goes to zero, these dynamics become those of TASEP.

Another feature one might include in a more realistic traffic model is the cascade effect of braking. The $q$-pushASEP includes this (Figure 10). Particles still jump right according to $q$-TASEP rules, however now particles may also jump left after exponential rate $L$ waiting times. When such a jump occurs, it prompts the next particle to the left to likewise jump left, with a probability given by $q^{\mathrm{gap}}$ where gap is the number of empty spaces between the original particle and its left neighbor. If that jump occurs, it may likewise prompt the next left particle to jump, and so on. Of course, braking is not the same as jumping backwards, however if one goes into a moving frame, this left jump is like a deceleration. It turns out that both of these models are solvable via the methods of Macdonald processes as well as stochastic quantum integrable systems and thusly it was been proved that, just as for ASEP, they demonstrate KPZ class fluctuation behavior (see the review [7, 15]).

## 7. Paths in a random environment

There is yet another class of probabilistic systems related to the corner growth model. Consider the totally asymmetric version of this model, started from wedge initial data. An alternative way to track the evolving height function is to record the time when a given box is grown. Using the labeling shown in Figure 11, let us call $L(x,y)$ this time, for $x,y$ positive integers. A box $(x,y)$ may grow, once its parent blocks $(x-1,y)$ and $(x,y-1)$ have both grown – though even then it must wait for an independent exponential waiting time which we denote by $w_{x,y}$. Thus $L(x,y)$ satisfies the recursion

$$L(x,y) = \max\bigl(L(x-1,y), L(x,y-1)\bigr) + w_{x,y}$$

subject to boundary conditions $L(x,0) \equiv 0$ and $L(0,y) \equiv 0$. Iterating yields

$$L(x,y) = \max_{\pi} \sum_{(i,j)\in\pi} w_{i,j}$$

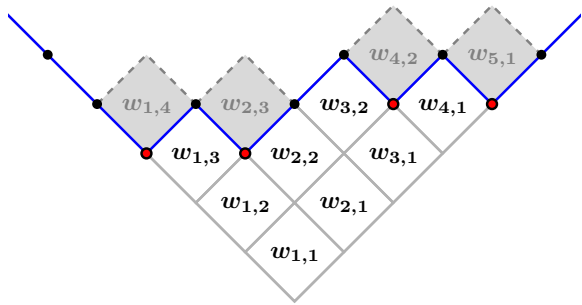

FIGURE 11. The relation between the corner growth model and last passage percolation with exponential weights. The $w_{i,j}$ are the waiting times between when a box can grow and when it does grow. $L(x,y)$ is the time when box $(x,y)$ grows.

where the maximum is over all up-right and up-left lattice paths between box $(1,1)$ and $(x,y)$. This model is called *last passage percolation* with exponential weights. Following from the earlier corner growth model results, one readily sees that for any positive real $(x,y)$, for large $t$, $L(\lfloor xt \rfloor, \lfloor yt \rfloor)$ demonstrated KPZ class fluctuations. A very compelling, and entirely open problem is to show that this type of behavior persists when the distribution of the $w_{i,j}$ is no longer exponential. The only other solvable case is that of geometric weights. A certain limit of the geometric weights leads to maximizing the number of Poisson points along directed paths. Fixing the total number of points, this becomes equivalent to finding the longest increasing subsequence of a random permutation. The KPZ class behavior for this version of last passage percolation was shown by Baik-Deift-Johansson [2].

There is another related integrable model which can be thought of as describing the optimal way to cross a large grid with stop lights are intersections. Consider the first quadrant of $\mathbb{Z}^2$ and to every vertex $(x,y)$ assign waiting times to the edges leaving the vertex rightwards and upwards. With probability $1/2$ the rightward edge has waiting time zero, while the upward edge has waiting time given by an exponential rate 1 random variables; otherwise reverse the situation. The edge waiting time represents the time needed to cross an intersection in the given direction (the walking time between lights has been subtracted). The minimal passage time from $(1,1)$ to $(x,y)$ is given by

$$P(x,y) = \min_\pi \sum_{e\in \pi} w_e$$

where $\pi$ goes right or up in each step and ends on the vertical line above $(x,y)$ and $w_e$ is the waiting time for edge $e\in \pi$. From the origin there will always be a path of zero waiting time, whose spatial distribution is that of the graph of a simple symmetric random walk. Just following this path one can get very close to the diagonal $x=y$ without waiting. On the other hand, for $x\neq y$, getting to $(\lfloor xt \rfloor, \lfloor yt \rfloor)$ for large $t$ requires some amount of waiting. Barraquand-Corwin [4] demonstrated that as long as $x\neq y$, $P(\lfloor xt \rfloor, \lfloor yt \rfloor)$ demonstrates KPZ class fluctuations. This should be true when $\pi$ is restricted to hit exactly $(x,y)$, though that result has not yet been proved. Achieving this optimal passage time requires some level of omnipotence as you must be able to look forward before choosing your route. As such, it could be considered as a benchmark against which to test various routing algorithms.

In addition to maximizing or minimizing path problems, the KPZ universality class describes fluctuations of 'positive temperature' version of these models in energetic or probabilistic favoritism is assigned to paths based on the sum of space-time random weights along its graph. One such system is called *directed polymers in random environment* and is the detropicalization of LPP where in the definition of $L(x,y)$ one replaces the operations of $(\max,+)$ by $(+,\times)$. Then the resulting (random) quantity is called the *partition function* for the model and its logarithm (the *free energy*) is conjectured for very general distributions on $w_{i,j}$ to show KPZ class fluctuations. There is one known integrable example of weights for which this has been proved – the inverse-gamma distribution, introduced by Seppäläinen [30] and proved in the work by Corwin-O'Connell-Seppäläinen-Zygouras [16] and Borodin-Corwin-Remenik [9].

The stop light system discussed above also has a positive temperature lifting of which we will describe a special case (see Figure 12 for an illustration). For each space-time vertex $(y,s)$ choose a random variable $u_{y,s}$ distributed uniformly on the interval $[0,1]$. Consider a random walk $X(t)$ which starts at $(0,0)$. If the random walk is in position $y$ at time $s$, then it jumps to position $y - 1$ at time $s + 1$ with probability $u_{y,s}$ and to position $y + 1$ with probability $1 - u_{y,s}$. With respect to the same environment of $u$'s, consider $N$ such random walks. The fact that the environment is fixed causes them to follow certain high probability channels. This type of system is called a *random walk in a space-time random environment* and the behavior of a single random walker is quite well understood. Let us, instead, consider the maximum of $N$ walkers in the same environment $M(t, N) = \max_{i=1}^{N} X^{(i)}(t)$. For a given environment, it is expected that $M(t, N)$ will localize near a given random environment dependent value. However, as the random environment varies, this localization value does as well in such a way that for $r \in (0, 1)$, and large $t$, $M(t, e^{rt})$ displays KPZ class fluctuations.

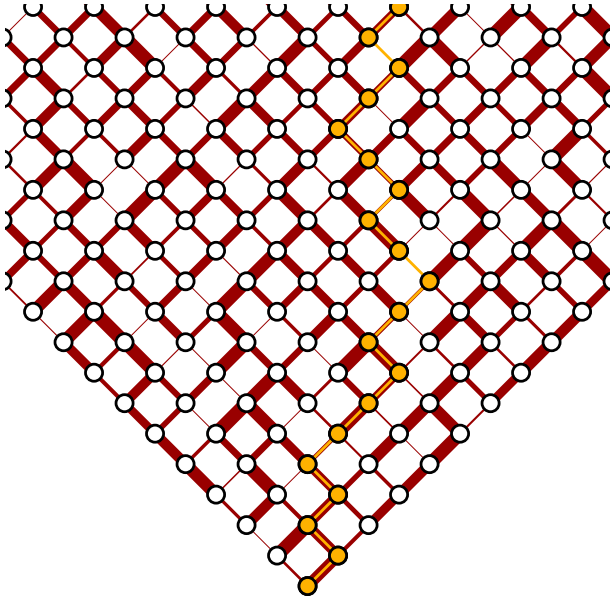

FIGURE 12. The random walk in a space time random environment. For each pair of up-left and up-right pointing edges leaving a vertex $(y, s)$, the width of the red edges is given by $u_{y,s}$ and $1 - u_{y,s}$ where $u_{y,s}$ are independent uniform random variables on the interval $[0, 1]$. A walker (the yellow highlighted path) then performs a random walk in this environment, jumping up-left or up-right from a vertex with probability equal to the width of the red edges.

## 8. BIG PROBLEMS

It took almost 200 years from the discovery of the Gaussian distributions to the first proof of its universality (the central limit theorem). So far, KPZ universality has withstood proof for almost three decades and shows no signs of yielding.

Besides universality, there remain a number of other big problems for which little to no progress has been made. All of the systems and results discussed herein have been $1 + 1$ dimensional, meaning that there is one time dimension and one space dimension. In the context of random growth, it makes perfect sense (and is quite important) to study surface growth $1 + 2$ dimensional. In the isotropic case (where the underly growth mechanism is roughly symmetric with respect to the two spatial dimensions) there are effectively no mathematical results though numerical simulations suggest that the $1/3$ exponent in the $t^{1/3}$ scaling for corner growth should be replaced by an exponent of roughly .24. In the anisotropic case there have been a few integrable examples discovered which suggest very different (logarithmic scale) fluctuations such as observed by Borodin-Ferrari [10].

Finally, despite the tremendous success in employing methods of integrable probability to expand and refine the KPZ universality class, there seems to still be quite a lot of room to grow and new integrable structures to employ. Within the physics literature, there are a number of exciting new directions in which the KPZ class has been pushed, including: out-of-equilibrium transform and energy transport with multiple conservation laws, front propagation equations, quantum localization with directed paths, and biostatistics. Equally important is to understand what type of perturbations break out of the KPZ class.

Given all of the rich mathematical predictions, one might hope that experiments have revealed the KPZ class behavior in nature. This is quite a challenge since determining scaling exponents and limiting fluctuations require immense numbers of repetition to experiments. However, there have been a few startling experimental confirmations of these behaviors in the context of liquid crystal growth, bacterial colony growth, coffee stains, and fire propagation (see [32] and references therein). Truly, the study of the KPZ universality class demonstrates the unity of mathematics and physics at its best.

8.0.1. *Acknowledgements.* The author appreciates the invitation from the AMS Notices editor Frank Morgan to write this article, as well as the collaborative work on figures with William Casselman. The author appreciates comments on a draft of this article by A. Borodin, P. Deift, P. Ferrari, and H. Spohn. This text is loosely based off his "Mathematic Park" lecture entitled *Universal phenomena in random systems* delivered at the Institut Henri Poincaré in May, 2015. He is partially supported by the NSF grant DMS-1208998, by a Clay Research Fellowship, by the Poincaré Chair, and by a Packard Fellowship for Science and Engineering.

I. Corwin, Columbia University, Department of Mathematics, 2990 Broadway, New York, NY 10027, USA, and Clay Mathematics Institute, 10 Memorial Blvd. Suite 902, Providence, RI 02903, USA, and Institut Henri Poincaré, 11 Rue Pierre et Marie Curie, 75005 Paris, France.

*E-mail address*: `ivan.corwin@gmail.com`